\documentclass{article}

\usepackage{amssymb}
\usepackage{makeidx}
\usepackage{multicol}
\usepackage{amsmath}
\usepackage[bottom]{footmisc}
\usepackage{graphicx}
\usepackage{color}
\usepackage{epstopdf}
\usepackage{url}

\newcommand{\al}[1]{\begin{align} #1 \end{align}}

\newcommand{\bsea}{\begin{subeqnarray}}
\newcommand{\esea}{\end{subeqnarray}}
\newcommand{\nn}{\nonumber}

\newcommand{\tr}{\mathop{\rm tr}}  

\newcommand{\Bc}{ \mathcal{B}}
\newcommand{\Cc}{ \mathcal{C}}
\newcommand{\Dc}{ \mathcal{D}}

\newcommand{\Gc}{ \mathcal{G}}
\newcommand{\Hc}{ \mathcal{H}}

\newcommand{\Lc}{ \mathcal{L}}

\newcommand{\Nc}{ \mathcal{N}}

\newcommand{\Qc}{ \mathcal{Q}}

\newcommand{\Sc}{ \mathcal{S}}

\newcommand{\Ds}{ \mathbb{D}}
\newcommand{\Es}{ \mathbb{E}}

\newcommand{\Ns}{ \mathbb{N}}

\newcommand{\Rs}{ \mathbb{R}}

\newcommand{\Ts}{ \mathbb{T}}

\newcommand{\Zs}{ \mathbb{Z}}

\def\qed{\hfill \vrule height 7pt width 7pt depth 0pt \smallskip}

\newcounter{pippo}

\newtheorem{defnn}{Definition}[section]
\newtheorem{remarkk}{Remark}[section]
\newtheorem{teor}{Theorem}[section]
\newtheorem{corr}{Corollary}[section]
\newtheorem{propo}{Proposition}[section]

\newtheorem{probl}[pippo]{Prfoblem}

\newcommand{\teo}{\begin{teor}}
\newcommand{\eteo}{\end{teor}}
\newcommand{\pb}{\begin{probl}}
\newcommand{\epb}{\end{probl}}
\newcommand{\df}{\begin{defnn}}
\newcommand{\edf}{\end{defnn}}
\newcommand{\aprop}{\begin{apropo}}
\newcommand{\eaprop}{\end{apropo}}
\newcommand{\alem}{\begin{alemm}}
\newcommand{\ealem}{\end{alemm}}

\begin{document}

\title{On the Robustness of the Bayes and Wiener Estimators under Model Uncertainty
\thanks{This work has been partially supported by the FIRB project ``Learning
meets time'' (RBFR12M3AC) funded by MIUR.}}

\author{Mattia Zorzi
\thanks{Dipartimento di Ingegneria dell'Informazione, Universit\`a degli studi di
Padova, via Gradenigo 6/B, 35131 Padova, Italy. Email: zorzimat@dei.unipd.it }}
\date{}
\maketitle

\begin{abstract}
This paper deals with the robust estimation problem of a signal given noisy observations. 
We assume that the actual statistics of the signal and observations belong to a ball about the nominal statistics.
This ball is formed by placing a bound on a suitable divergence (or distance) between the actual and the nominal statistics.   
Then, the robust estimator is obtained by minimizing the mean square error according to the least favorable statistics 
in that ball. Therefore, we obtain a divergence-based minimax approach to robust estimation. Choosing a set of divergences, called Tau divergence family,
 we show that the Bayes estimator based on the nominal statistics is the optimal solution. 
Moreover, in the dynamic case, the optimal offline estimator is the noncausal Wiener filter based on the nominal statistics.

Keywords: Robust filtering, minimax problem, Tau divergence family, risk-sensitive estimation problem, minimum entropy problem. 
 \end{abstract}

\section{Introduction}\label{sec_introduction}

Consider the problem of estimating a signal from noisy observations. Typically the actual statistics of the signal and observations are known only imprecisely that is only the nominal statistics are known. In this situation, one would wonder 
how to construct an estimator which is robust to this model uncertainty. According to the minimax approach, the statistics of the signal and observations are assumed to belong to a neighborhood of the nominal ones. Then, the robust estimator is characterized by a minimax problem consisting of finding the estimator which minimizes the mean square error for the least favorable statistics in the neighborhood. The latter can be specified in different ways, e.g. it can be based on a $\epsilon$-contamination model, a total variation model or a spectral band model \cite{KASSAM_POOR_1985}. \cite{ROBUSTNESS_HANSENSARGENT_2008} propose to specify this neighborhood through an
uncertainty ball which is formed by placing a bound (i.e. tolerance) $c$ on the {\em Kullback-Leibler}  divergence between the actual and the nominal statistics. This characterization  of the uncertainty is supported by the fact that nominal models are identified from data according to the maximum-likelihood principle which turns out to be equivalent to the minimization of the {\em Kullback-Leibler}  divergence over a suitable parametric model class. Accordingly, the maximum-likelihood approach provides the nominal model and the bound $c$. \cite{LEVY_NIKOUKHAH_2004} showed that if the nominal statistics is Gaussian then the least favorable statistics is Gaussian and the robust estimator coincides with the nominal {\em Bayes}  estimator, i.e. the {\em Bayes} estimator based on the nominal statistics.

In this paper we consider a family of uncertainty balls which are formed by placing a bound on a set of divergences, called $\tau$ divergence family, between Gaussian statistics. Each divergence of this set is characterized by parameter $\tau$ which is a real number belonging to the interval $[0,1]$. In particular, for $\tau=0$ we obtain the {\em Kullback-Leibler} divergence.
This characterization  of the uncertainty is supported by the fact that recently it has been proposed a system identification procedure which finds the model from data by minimizing the $\tau$ divergence family over a suitable parametric model class \cite{DUAL}. Accordingly, this system identification procedure provides the nominal model and the bound $c$ for the proposed uncertainty ball. It turns out that the (family of) robust estimators, solution to the minimax problem with this family of uncertainty balls, coincide with the {\em Bayes} estimator
based on the nominal statistics. Accordingly, the nominal {\em Bayes} estimator is robust in this wide family of (Gaussian) uncertainty classes. Our result also gives the analytical form of the least favorable statistics of the estimation error. 
This result allows to define a new family of robust Kalman filters obtained by iterating the {\em Bayes} estimator with 
the least favorable statistics \cite{STATETAU,ROBUST_STATE_SPACE_LEVY_NIKOUKHAH_2013}.
Our minimax approach can be also relaxed expressing the $\tau$ divergence constraint as a soft one. We will show that this relaxation corresponds to a new family of risk-sensitive problems (in the sense of \cite{boel2002robustness}) which is also linked to a new family of minimum entropy problems (in the sense of \cite{MUSTAFA_GLOVER_1990}). These results can be extended in the dynamic case: We will prove that the nominal noncausal {\em Wiener} filter is the optimal robust offline estimator according to the minimax approach based on the $\tau$ divergence family. Finally, through a simulation study we analyze the features of the least favorable statistics of the estimation error for these uncertainty classes. Simulations show that the parameter 
$\tau$ tunes how the uncertainty is spread among the components (components and frequencies for the dynamic case) of the least favorable statistics.

The outline of the paper is as follows. In Section \ref{sec_tau_div} we define  the $\tau$ divergence family for Gaussian random vectors and stationary stochastic processes.  Section \ref{sec_static} deals with the static case: We show that the nominal {\em Bayes} estimator is the solution to the family of divergence-based minimax approaches.  Section \ref{sec_wiener} deals with the dynamic case: We show that the nominal noncausal {\em Wiener} estimator is the solution to the family of divergence-based minimax approaches. In Section \ref{sec_sim}, we present the simulation study.
Finally, in Section \ref{sec_concl} we draw the conclusions. In order to streamline the presentation all the proofs are deferred to the Appendix.

In the paper, we will use the following conventions. $\Zs,\Ns$ and $\Rs$ denote the set of integer, natural and real numbers, respectively. $\Rs_+$ denotes the set of positive real numbers.
Given $x\in\Rs^q$, $\|x\|$ denotes its {\em Euclidean} vector norm. Moreover, $\|x\|_K=\sqrt{x^T Kx}$. $I_q$ denotes the identity matrix of dimension $q$.  $\tr(P)$ denotes the trace of matrix $P$. $\Qc^{q}$ denotes
the vector space of symmetric matrices of
dimension $q\times q$. The $i$-th singular value of $P\in\Qc^q$ is denoted by $\sigma_i(P)$ and we assume that
$\sigma_1(P)\geq \sigma_2(P)\geq \ldots \sigma_q(P)$. $\|P\|$ denotes the spectral norm of $P$, that is $\|P\|=\sigma_1(P)$. $\Qc_+^q$ is
the cone of positive definite symmetric matrices of dimension $q\times q$. $\log(P)$ denotes the logarithm of matrix $P$ and $P^\tau$ is the $\tau$-th power of matrix $P$. Matrix
functions defined over the unit circle $\Ts=\{e^{j\vartheta}\, : \, \vartheta\in[0,2\pi]\}$ are denoted by capital Greek
letters and the dependence upon $\vartheta$ is sometime dropped to simplify the notation, i.e. $\Sigma$ instead of
$\Sigma(\vartheta)$. A star denotes transposition plus conjugation, that is $\Sigma(\vartheta)^*=\Sigma(-\vartheta)^T$.
$\|\Sigma\|_\infty$ is the infinity matrix norm
of $\Sigma$, that is $\|\Sigma\|_\infty=\sup_\vartheta \,\sigma_1(\Sigma(\vartheta))$. 
$\Qc^{q}(\Ts)$ denotes the vector space of para-symmetric matrix functions of
dimension $q\times q$, i.e. if $\Sigma\in\Qc^{q}(\Ts)$ then $\Sigma=\Sigma^*$. $\Qc_{+}^q(\Ts)$ denotes the cone of para-symmetric matrix functions of dimension $q\times q$ which are positive definite over $\Ts$. Given
$\Sigma\in \Qc^q({\Ts})$, the shorthand notation  $\int \Sigma$ means the integration over the unit circle with respect to the
normalized {\em Lebesgue} measure, that is $(2\pi)^{-1}\int_0^{2\pi} \Sigma(\vartheta)\mathrm{d}\vartheta$. $\delta(\vartheta)$ denotes the {\em Dirac delta function}.

\section{$\tau$ Divergence for Gaussian Vectors and Processes}\label{sec_tau_div}
Let $z$ be a Gaussian random vector of dimension $q$ with probability density 
\al{\label{ftilda}\tilde f(z)&=\nn\\ & \frac{1}{\sqrt{(2\pi)^{q} \det \tilde K_z }}\exp\left(-\frac{1}{2}(z-\tilde m_z)^T \tilde K_z^{-1} (z-\tilde m_z)\right),} 
where $\tilde m_z\in\Rs^{q}$ and
$\tilde K_z\in\Qc^{q}_+$.  Let $\hat z$ denote the minimum variance predictor of $z$ based on the nominal probability density $f$ 
\al{\label{fnomi} f(z)& =\nn\\ &\frac{1}{\sqrt{(2\pi)^{q} \det K_z }}\exp\left(-\frac{1}{2}(z-m_z)^T K_z^{-1} (z-m_z)\right) }
with $m_z\in\Rs^{q}$ and $K_z\in\Qc^{q}_+$. Thus, $\hat z=m_z$. Let $e=z-\hat z$ be the corresponding innovation vector.     
Accordingly, $e^N=L_z^{-1}(z-\hat z)$ is the normalized innovation vector with $L_z$ a square root of $K_z$, i.e. $K_z=L_zL_z^T$. It is not difficult to see that $e^N$ is Gaussian with mean $m=L_z^{-1}\Delta m_z$, $\Delta m_z=\tilde m_z-m_z$, and covariance matrix $K=L_z^{-1} \tilde{K}_z L_z^{-T}$. If $\tilde f$ coincides with $f$ then we have that $e^N$ is with zero mean and identically independently distributed components, i.e. 
$m=0$ and $K=I$. Therefore, $(m,K)$ represents a mismatch criterium which naturally occurs in prediction error estimation \cite{LINDQUIST_PICCI}. This leads us to measure the mismatch between $\tilde f$ and $f$ by quantifying the mismatch between $(m,K)$ and $(0,I)$:
\al{\label{D_pred}\Dc(\tilde f \| f)= \ell(m,K),}
where $\ell: \Rs^{q}  \times \Qc_+^q \rightarrow \Rs \cup \{\infty\}$ is a function such that $\ell\geq 0$   and equality holds if and only if $m=0$ and $K=I$.
We consider the following function parametrized by $\tau\in[0,1]$:
\al{	\label{ell_tau} &\ell(m,K)=\nn\\
& \left\{
                               \begin{array}{ll}
                                 \| m\|^2+\tr\left(-\log(K)+ K-I_{q} \right), & \hbox{$\tau=0$} \\
 \frac{1}{1-\tau}\| m\|^2+\tr \left(\frac{1}{\tau(\tau-1)} K^\tau+\frac{1}{1-\tau} K+\frac{1}{\tau} I_q\right), & \hspace{-0.1cm}\hbox{$0<\tau<1$} \\
                               \delta( m)+  \tr\left( K \log(K)-K+I_{q}\right), & \hbox{$\tau=1$.}
                               \end{array}
                             \right. }
 There are several ways to construct $\ell$. We will motivate our choice in Remark \ref{rem_scelta}.
Substituting (\ref{ell_tau}) in (\ref{D_pred}) we obtain the following family of divergences indexed by $\tau$: 
\al{ \label{def_tau_div}& \Dc_\tau(\tilde f\| f)=\nn\\
& \left\{
                               \begin{array}{ll}
                                 \|\Delta m_z\|^2_{K_Z^{-1}}+\tr\left(-\log(\tilde K_z K_z^{-1})\right. &\\ \hspace{0.3cm}\left.+\tilde K_z K_z^{-1}-I_{q} \right), & \hbox{$\tau=0$} \\
                                	\frac{1}{1-\tau} \|\Delta m_z\|^2_{K_Z^{-1}}+\tr\left(\frac{1}{\tau(\tau-1)}(L_z^{-1}\tilde K_z L_z^{-T})^{\tau}\right. & \\ \hspace{0.3cm}\left.  +\frac{1}{1-\tau} \tilde K_z K_z^{-1}+\frac{1}{\tau}I_{q}\right), & \hspace{-0.4cm}\hbox{$0<\tau<1$} \\
                               \delta(\Delta m_z)+  \tr\left( L_z^{-1}\tilde K_z L_z^{-T}\log(L_z^{-1}\tilde K_z                                    L_z^{-T})\right. &\\
                               \hspace{0.3cm}\left.-\tilde K_z K_z^{-1}+I_{q}\right), & \hbox{$\tau=1$. }
                               \end{array}
                             \right. }
It is worth noting that (\ref{def_tau_div}) coincides with the $\tau$ divergence between covariance matrices \cite{OPTIMAL_PREDICTION_ZORZI_2014} when $\tilde f$ and $f$ have the same mean. Moreover, (\ref{def_tau_div}) coincides with the {\em Kullback-Leibler} divergence for $\tau=0$.

\begin{propo} \label{prop_tau_divergence}$\Dc_\tau(\tilde f \| f)\geq 0$ and equality holds if and only if $\tilde f=f$. \end{propo}

The mismatch criterium above can be extended to the dynamic case. Let $z(t)$ be a stationary Gaussian process defined over $t\in \Zs$ of dimension $q$
with probability measure $\tilde f$. The latter is completely characterized  by its power spectral density 
$\tilde S_z(\vartheta)=2\pi \tilde m_z \tilde m_z^T\delta(\vartheta)+ \tilde \Sigma_z(\vartheta)$,
where $\tilde m_z\in\Rs^q$ is the mean and $\tilde \Sigma_z\in\Qc_+^q(\Ts)$ is the discrete time {\em Fourier} transform of the covariance matrix function 
\al{ K_z(s)&=\Es_{\tilde f}[(z(t)-\tilde m_z)(z(t-s)-\tilde m_z)^T],\; s\in\Zs.\nn} 
Let $\hat z(t)$ be the minimum variance linear one-step-ahead predictor based on the nominal probability measure $f$ with power spectral density $S_z(\vartheta)=2\pi m_z m_z^T\delta(\vartheta)+\Sigma_z(\vartheta)$.
It is not difficult to see that corresponding normalized innovation process $e^N(t)$ is stationary Gaussian with power spectral density 
$S(\vartheta)=m m^T\delta(\vartheta)+\Sigma(\vartheta)$, where $\Sigma_z=\Gamma_z^{-1} \tilde{\Sigma}_z\Gamma_z^{-*}$, $\Gamma_z$ is a left squared spectral factor of $\Sigma_z$, i.e. $\Sigma_z=\Gamma_z\Gamma_z^*$, $m=\Gamma_z(0)^{-1}\Delta m_z$, $\Delta m_z=\tilde{m}_z-m_z$. Clearly, the more $e^N(t)$ is similar to white Gaussian noise, i.e. $S\approx I$,  the closer $\tilde f$ and $f$ are. Accordingly, we measure their mismatch as follows 
$\Sc(\tilde f \| f)=\int  \ell(2\pi m,\Sigma).$
Choosing $\ell$ as in (\ref{ell_tau}), we obtain the following family of divergences indexed by $\tau$: 
\al{ \label{def_tau_div_spectrum} &\Sc_\tau(\tilde f\| f)=\nn\\
&\hspace{0.1cm} \left\{
                               \begin{array}{ll}
                                 \|\Delta m_z\|_{\Sigma_z(0)^{-1}}^2+\int \tr \left(-\log(\tilde \Sigma_z \Sigma_z^{-1})\right. & \\  \hspace{0.3cm}\left.+\tilde \Sigma_z \Sigma_z^{-1}-I_{q} \right), & \hspace{-0.2cm}\hbox{$\tau=0$} \\
                                 \frac{1}{1-\tau} \|\Delta m_z\|_{\Sigma_z(0)^{-1}}^2+\int \tr \left(\frac{1}{\tau(\tau-1)}(\Gamma_z^{-1}\tilde \Sigma_z \Gamma_z^{-*})^{\tau}\right. & \\ \hspace{0.3cm} \left.
                                 +\frac{1}{1-\tau}\tilde \Sigma_z \Sigma_z^{-1}+\frac{1}{\tau}I_{q}\right), &\hspace{-0.8cm} \hbox{$0<\tau<1$} \\
                                \delta(\Delta m_z)+\int  \tr\left( \Gamma_z^{-1}\tilde \Sigma_z \Gamma_z^{-*}\log(\Gamma_z^{-1}\tilde \Sigma_z \Gamma_z^{-*}) \right. & \\ \hspace{0.3cm}\left. -\tilde \Sigma_z \Sigma_z^{-1}+I_{q}\right), &\hspace{-0.3cm} \hbox{$\tau=1$.}
                               \end{array}
                             \right. } Note that, for $\tau=0$ we obtain the {\em Itakura-Saito} distance, \cite{itakura1968analysis}. In the case that $\tilde f$ and $f$ have the same mean, we obtain the $\tau$ divergence defined in \cite{OPTIMAL_PREDICTION_ZORZI_2014}.
\begin{propo} \label{prop_d_tau_dyn}$\Sc_\tau(\tilde f\| f)\geq 0$ and equality holds if and only if $\tilde f=f$. \end{propo}

\section{Robust Static Estimation} \label{sec_static}
We consider a static estimation problem where we seek to estimate a random vector $x\in\Rs^{n}$ given an observation $y\in\Rs^{p}$.
Assume the joint vector $z:= [ x^T; y^T]^T$ is Gaussian with nominal probability density $f$ defined in (\ref{fnomi})
where $q=n+p$.
We conformably partition the mean vector and the covariance matrix of $z$ according to $x$
and $y$:
\al{ m_z=\left[
          \begin{array}{c}
            m_x \\
            m_y \\
          \end{array}
        \right], \;\; K_z=\left[
                            \begin{array}{cc}
                              K_x & K_{xy} \\
                              K_{yx} & K_y\\
                            \end{array}
                          \right]
\nn.} Let $\tilde f(z)$ be the actual probability density defined in (\ref{ftilda})
where $q=n+p$.
We consider  the closed ball centered on $f$:
\al{ \label{def_Bc_static}\Bc_\tau:=\{ \tilde f   \hbox{ s.t. } \Dc_\tau(\tilde f\| f)\leq c\},}
where $\Dc_\tau$ has been defined in (\ref{def_tau_div}), $c\in\Rs_+$ is a fixed tolerance which accounts for the maximum allowable deviance. 
Therefore, the hope is that $\Bc_\tau$ contains the actual (unknown) probability density $\tilde f$.
 Note that, $\Bc_\tau$ depends on $\tau\in[0,1]$. Accordingly, by changing $\tau$ the set of all possible probability densities changes.
In this way, we have a family of uncertainty classes parametrized by $\tau$. It is worth noting that $\mathcal{B}_\tau$ also depends on $c$. However, we dropped this dependence to ease the notation.  

\begin{remarkk} The nominal probability density function can be identified from data by solving $f=\mathrm{argmin}_{\bar f\in \Cc} \Dc_\tau( f_{S}\|\bar f)$ where $\Cc$ is a suitable parametric family of probability densities  and $f_{S}$ is Gaussian with mean and covariance matrix, respectively, the sample mean and the sample covariance matrix computed from the data \cite{DUAL}. Since $f_S$ is the best probability density fitting the data, it is then realistic to specify the uncertainty with $\Bc_{\tau}$ with $c=\Dc_\tau(f_{S}\| f)$.
\end{remarkk}

We shall use the minimax viewpoint to design our robust estimator of $x$ \cite{LEVY_NIKOUKHAH_2004,ROBUSTNESS_HANSENSARGENT_2008}. More precisely, whenever
we seek to design an estimator minimizing a suitable loss function, an hostile player, say ``nature'', conspires to select the worst possible probability density in $\Bc_\tau$.
Let $g(y)$ denote an estimator of $x$ based on the observation $y$. We evaluate its performance  through the mean square error \al{  J(\tilde f,g)&= \Es_{\tilde f}[\| x-g(y)\|^2]=\int _{\Rs^{n+p}} \| x-g(y)\|^2 \tilde f(z)\mathrm{d} z. \nn }
Let $\Gc$ denote the set of estimators $g(y)$ such that $\Es_{\tilde f}[\|g(y)\|^2]$ is finite for any $\tilde f \in\Bc_\tau$.
Our optimal robust estimator is the solution to the following minimax problem
 \al{ \label{minimax_problem_fg}\underset{g\in\Gc}{\min}\;\underset{\tilde f\in\Bc_\tau}{\max} \;J(\tilde f,g).}

\teo \label{theorem_static} The optimal robust estimator, according to (\ref{minimax_problem_fg}), is the {\em Bayes} estimator based on $f $ \al{ \label{Opt_estimator}g^\circ(y)=G^\circ(y-m_y)+m_x}
with $ G^\circ=K_{xy}K_y^{-1}$.
The least favorable probability density $\tilde f^\circ$ has mean vector and covariance matrix  \al{ \label{def_K_tilda} \tilde m_z^\circ =m_z,\;\;\tilde K^\circ_z=\left[
                                                                         \begin{array}{cc}
                                                                           \tilde K_x & K_{xy} \\
                                                                           K_{yx} &  K_y \\
                                                                         \end{array}
                                                                       \right],
} wherein only the covariance of $x$ is perturbed with respect to the nominal covariance matrix. The nominal and the least favorable estimation error have zero mean and covariance matrix, respectively,  \al{ P =  K_x-K_{xy}K_y^{-1}K_{yx} ,\;\;    \tilde P= \tilde K_x-K_{xy}K_y^{-1}K_{yx}.\nn}
Moreover,
\al{  \label{relaz_delle_P}\tilde P=\left\{
                                   \begin{array}{ll}
                                     L_P\left(I_{n}-\frac{1-\tau}{\lambda}L_P^TL_P\right)^{\frac{1}{\tau-1}}L_P^T, & 0 \leq \tau<1 \\
                  L_P\exp\left(\frac{1}{\lambda} L_P^TL_P\right) L_P^T                   , & \tau=1,
                                   \end{array}
                                 \right. }
where $P=L_PL_P^T$ and $\lambda$, with $\lambda>(1-\tau)\|P\|$, is the unique {\em Lagrange} multiplier
 such that $\Dc_\tau(\tilde f\| f )=c$. \eteo

Theorem \ref{theorem_static} shows that the {\em Bayes} estimator based on the nominal statistics $f$ is robust with respect to the $\tau$ divergence constraint.
The worst situation occurs when all the mismodeling budget is allocated in a perturbation of the covariance matrix $K_x$. In Problem (\ref{minimax_problem_fg})
several divergence families can be used to characterize $\Bc_\tau$, such as the $\alpha$ divergence, \cite{ALPHA}, and the $\beta$ divergence family \cite{cichocki2010families}. Although the existence of the solution to (\ref{minimax_problem_fg}) with those uncertainty classes is guaranteed, the optimal Bayes estimator is not necessarily based on the nominal statistic and such a solution does not admit a closed form. The mean square error (MSE) corresponding to the nominal probability density is
$ \mathrm{MSE}=\Es_f[\|e\|^2]=\tr(P)$, while the MSE corresponding to the least favorable probability density $\tilde f^\circ$
is $ \widetilde{\mathrm{MSE}}= \Es_{\tilde f^\circ}[\|e\|^2]=\tr(\tilde P)$.
In view of (\ref{relaz_delle_P}), it follows that $\tilde P- P\in\Qc^n_+$ therefore $ \widetilde{\mathrm{MSE}}> \mathrm{MSE}$ and the
additional MSE occasioned by the least favorable model perturbation is
$ \Delta\mathrm{MSE}=\tr(\tilde P-P).$

\begin{remarkk} \label{rem_scelta}One would wonder why in (\ref{def_tau_div}) we consider $\frac{1}{1-\tau}\|\Delta m_z\|^2_{K_z^{-1}}$ instead of the
simpler term $\|\Delta m_z\|^2_{K_z^{-1}}$. Indeed, it is not difficult to see that Theorem \ref{theorem_static} still holds with
$\frac{1}{1-\tau}\|\Delta m_z\|^2_{K_z^{-1}}$ replaced by $\|\Delta m_z\|^2_{K_z^{-1}}$ in (\ref{def_tau_div}).
The unique difference is that $\lambda$ must be such that $\lambda>\|P\|$ for any $\tau\in[0,1]$.
Therefore, with the choice
$\|\Delta m_z\|^2_{K_z^{-1}}$ we restrict the allowable values for the Lagrange multiplier $\lambda$ and thus the allowable
least favorable covariance matrices $\tilde K_z$.
Accordingly, taking the term $\| \Delta m_z\|_{K_z}^2$, the freedom of the nature (i.e. the hostile player) is restricted.
\end{remarkk}

It is worth comparing our result with the one in \cite{LEVY_NIKOUKHAH_2004}.  Theorem \ref{theorem_static} generalizes the case
$\tau=0$ analyzed in Theorem 1 in \cite{LEVY_NIKOUKHAH_2004}. On the other hand,  Theorem 1
shows the least square estimator is robust over the more general allowable set
$\{\tilde f \,:\, \Ds_{KL}(\tilde f\|f)\leq c\}$,
where $\Ds_{KL}(\tilde f\|f)$
is the {\em Kullback-Leibler} divergence among probability densities, and $\tilde f$ is not necessarily Gaussian. However, our result cannot be extended to such general case because $\Dc_\tau$
is a divergence family which only measures the deviation among Gaussian vectors.

\begin{figure}[htbp]
\begin{center}
\includegraphics[width=\columnwidth]{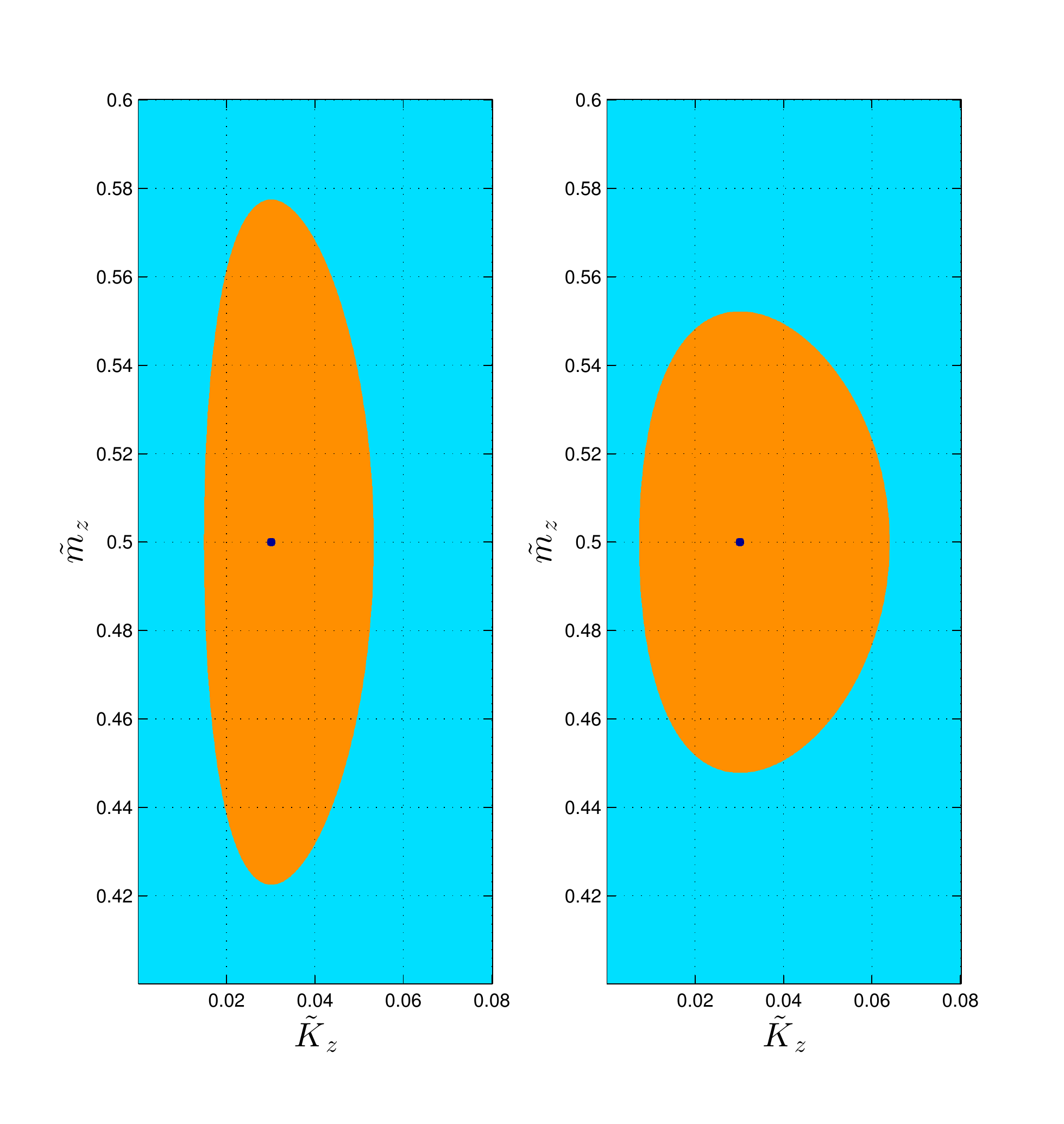}
\end{center}
 \caption{In orange the uncertainty ball with $\tau=0$, $c=0.2$ (left panel) and $\tau=0.8$ and $c=0.454$ (right panel); the blue point is the nominal statistics.}\label{levelsetfig}
\end{figure}

In oder to understand the influence of parameter $\tau$ on the uncertainty ball $\Bc_\tau$, we consider the case in which $z$ is a Gaussian random variable, i.e. $q=1$, with nominal mean $m_z=0.5$ and nominal variance $K_z=0.03$. We consider $\Bc_\tau$ with $\tau=0$ and $c=0$, left panel of Figure \ref{levelsetfig}, and $\Bc_\tau$ with $\tau=0.8$ and $c=0.454$, right panel of Figure \ref{levelsetfig}.    
 The tolerance $c$ for the two balls is chosen in such a way that their measure is the same. As we can see, increasing $\tau$ the uncertainty increases for the variance while it decreases for the mean.
This observation holds also with $q>1$, indeed the first term in (\ref{def_tau_div}) measures the deviance between the actual and the nominal mean, and, increasing $\tau$, this term becomes preponderant than the second one when $m_z\neq \tilde m_z$. In other words through parameter $\tau$ we tune how to allocate the mismodeling budget between the mean and the covariance matrix.

Problem (\ref{minimax_problem_fg}) can be relaxed in the following way 
 \al{ \label{minimax_problem_fg_relaxed}\underset{g\in\Gc}{\min}\;\underset{\tilde f\in\Bc_\tau^\infty}{\max} \;J(\tilde f,g)+\lambda(c-\Dc_\tau(\tilde f\| f)),}
where $ \Bc_\tau^\infty=\{ \tilde f \hbox{ s.t. } \Dc_\tau (\tilde f\| f)<\infty\} $
and $\Gc$ is the set of all estimators such that $\Es_{\tilde f}[\|g(y)\|^2]$ is finite for any $\tilde f\in\Bc_\tau^\infty$.
Here, $\lambda>0$ is {\em a priori} fixed and such that $\lambda>(1-\tau)\| P\|$. In this way the mismodeling  tolerance is expressed as a soft constraint adding the penalty term $\lambda(c-\Dc_\tau(\tilde f\| f))$.  
\begin{corr} The optimal estimator, according to (\ref{minimax_problem_fg_relaxed}), is still the nominal {\em Bayes} estimator. The least favorable probability density $\tilde f^\circ$ has mean vector $\tilde m_z=m_z$ and covariance matrix $\tilde K_z^\circ$ as in (\ref{def_K_tilda}). The least favorable estimation error has zero mean and covariance matrix $\tilde P$ as in (\ref{relaz_delle_P}) where $\lambda$ now has been chosen {\em a priori}. \end{corr}

In the perspective presented in \cite{boel2002robustness,HANSEN_SARGENT_2005,HANSEN_SARGENT_2007,ROBUSTNESS_HANSENSARGENT_2008,LEVY_NIKOUKHAH_2004}, Problem (\ref{minimax_problem_fg_relaxed})
represents a generalization of the risk-sensitive static estimation problem \cite{RISK_WHITTLE_1980,speyer2008stochastic,H_INF_HASSIBI_SAYED_KAILATH_1999,zhou1996robust,ZORZI_CONTRACTION_CDC,LEVY_ZORZI_RISK_CONTRACTION}. Thus, the nominal {\em Bayes} estimator is also 
optimal for the $\tau$ risk-sensitive static estimation problem.

Now we show that Problem (\ref{minimax_problem_fg_relaxed}) is the solution to a new minimum entropy problem in the sense of \cite{MUSTAFA_GLOVER_1990}. Let $g(y)=Gy+h$ be an estimator of $x$. Under the nominal model, the estimation error $e=x-g(y)$ is Gaussian with
\al{ m_e&=\Es_f[e]=\left[
                  \begin{array}{cc}
                    I_n & -G\\
                  \end{array}
                \right] m_z-h\nn\\
K_e&=\left[
                  \begin{array}{cc}
                    I_n & -G \\
                  \end{array}
                \right] K_z \left[
                                   \begin{array}{c}
                                     I_n \\
                                     -G^T \\
                                   \end{array}
                                 \right]. \nn }
Note that, the mean and covariance matrix of $e$ depends on $g(y)$. The idea is to characterize the robust estimator through the following minimum entropy problem
\al{ \label{min_entropy_pb}\underset{g\in\Gc}{\min} \;\Hc_\tau(e,\lambda), }
where $\Hc_\tau$ is an entropy-like function which guarantees that the mean and the covariance matrix of $e$ are bounded in some sense. Such boundedness is tuned by parameter $\lambda\in\Rs_+$. Next, we characterize $\Hc_\tau$.
\df The $\tau$ entropy
family of $e$ is defined as
\al{  & \Hc_\tau(e,\lambda)=\\
&\left\{\begin{array}{ll}  m_e^T \left( I_n -\frac{1}{\lambda}K_e\right)^{-1} m_e -\lambda\log\det \left(I_n-\frac{1}{\lambda}K_e\right),
& \tau=0 \\ 
m_e^T\left(I_n-\frac{1-\tau}{\lambda} K_e\right)^{-1}m_e \\ \hspace{0.3cm}    +\frac{\lambda}{\tau}\tr\left( \left(I_n-\frac{1-\tau}{\lambda}K_e\right)^{\frac{\tau}{\tau-1}}-I_n\right),
 & \hspace{-0.6cm} 0<\tau<1\\ 
m_e m_e^T+ \lambda\tr\left(\exp\left(\frac{1}{\lambda}K_e\right)-I_n\right), & \tau=1\end{array}\right. \nn }
for $\lambda>(1-\tau)\|K_e\|$ otherwise $\Hc_\tau(e,\lambda)=\infty$.  \edf
It is not difficult to see that $\Hc_\tau$ is continuous with respect to $\tau\in[0,1]$.
Note that, ${\Hc}_0$ is the usual entropy, see \cite{MUSTAFA_GLOVER_1990}, used in $H_\infty$ control.  

\begin{propo}  \label{prop_tau_entropy} The following properties hold:
\begin{itemize}
  \item $\Hc_\tau(e,\lambda)\geq 0$ and equality holds if and only if $e=0$
\item $\Hc_\tau(e,\cdot)$ is a monotone decreasing function.
\end{itemize}
\end{propo}
 
The Theorem below states the connection between our minimax approach for robust estimation and the minimum entropy estimation.
\teo \label{TH_min_entropy} Problem (\ref{minimax_problem_fg_relaxed}) and Problem (\ref{min_entropy_pb}) are equivalent. 
\eteo

\section{Noncausal Robust Filtering}\label{sec_wiener}
Let $x(t)$ and $y(t)$ be two jointly stationary Gaussian processes defined over $t\in\Zs$ of dimension
$n$ and $p$, respectively. We consider the noncausal filtering problem, that is to estimate
$x(t)$ given the observations $\{y(s),\; s\in\Zs\}$. We define
$z(t)=[ x(t)^T\;\; y(t)^T]^T$. The nominal probability measure of $z(t)$, say $f$, has power spectral density 
$S_z(\vartheta)=2\pi m_z m_z^T\delta (\vartheta)+ \Sigma_z(\vartheta)$, where
\al{ m_z=\left[\begin{array}{c} m_x \\ m_y\end{array}\right],\;\; \Sigma_z(\vartheta)=\left[
         \begin{array}{cc}
            \Sigma_x(\vartheta) & \Sigma_{xy}(\vartheta) \\
           \Sigma_{yx}(\vartheta) & \Sigma_y(\vartheta) \\
         \end{array}
       \right] .\nn}
The actual one, say $\tilde f$, has power spectral density $\tilde S_z(\vartheta)=2\pi \tilde m_z \tilde m_z^T\delta (\vartheta)+ \tilde \Sigma_z(\vartheta)$, where
\al{ \tilde m_z=\left[\begin{array}{c} \tilde m_x \\ \tilde m_y\end{array}\right],\;\; \tilde \Sigma_z(\vartheta)&=\left[
         \begin{array}{cc}
           \tilde \Sigma_x(\vartheta) & \tilde \Sigma_{xy}(\vartheta) \\
          \tilde  \Sigma_{yx} (\vartheta) & \tilde \Sigma_y(\vartheta) \\
         \end{array}
       \right].\nn}
 Suppose that the actual power spectral density belongs to the following closed ball parametrized by $\tau$
\al{ \label{def_Bc_dynamic}\Bc_\tau=\{\tilde f  \hbox{ s.t. } \Sc_\tau(\tilde f\| f)\leq c\},}
where $\Sc_\tau$ has been defined in (\ref{def_tau_div_spectrum}) and $c\in\Rs_+$ is a fixed tolerance. 
Similarly to the static case, we design the robust noncausal filter according to the minimax point of view.
Let $\Gc$ denote class of estimators of $x(t)$ having the following structure
\al{ g(y,t)=\sum_{k=-\infty}^\infty G_ky(t-k)+h,}
where the filter $\Lambda(\vartheta)=\sum_{k=-\infty}^\infty G_k e^{-j\vartheta k}$
is Bounded Input Bounded Output (BIBO) stable. Our robust noncausal filter is the solution to the following minimax problem
\al{ \label{minimax_wiener}\underset{g\in\Gc}{\min}\;\underset{\tilde f\in\Bc_\tau}{\max} \;J(\tilde f,g),}
where $J(\tilde f,g)= \Es_{\tilde f}[\|x(t)-g(y,t)\|^2]$.

\teo \label{thm_wiener_saddle_point} Let $\tau$ be such that $\frac{1}{1-\tau}\in\Ns$ and $\Sigma_z$ have bounded {\em McMillan} degree. The robust estimator $g^\circ(y,t)$, according to (\ref{minimax_wiener}), is the noncausal {\em Wiener} filter based on $f$, that is
\al{ h =\left[\begin{array}{cc} 1 &  -\Lambda  (0)\\ \end{array}\right]m_z,\;\; \Lambda =\Sigma_{xy}\Sigma_y^{-1}.\nn}
The least favorable probability measure $\tilde f^\circ$ has power spectral density $\tilde S_z^\circ(\vartheta)=2\pi m_z m_z^T\delta(\vartheta)+\tilde \Sigma_z^\circ(\vartheta) $
with \al{  \tilde  \Sigma_z^\circ =\left[        \begin{array}{cc}
                                         \tilde\Sigma_x  & \Sigma_{xy}  \\
                                         \Sigma_{yx}  & \Sigma_y  \\
                                       \end{array}
                                     \right], \nn
} where only $\tilde \Sigma_x $ is perturbed. The power spectral density of the estimation error with respect to the nominal
and the least favorable probability measure, respectively, are  
\al{ \Sigma_e =\Sigma_x-\Sigma_{xy}\Sigma_y^{-1}
\Sigma_{yx}, \;\; \tilde \Sigma_e=\tilde \Sigma_x-\Sigma_{xy}\Sigma_y^{-1}
\Sigma_{yx}.\nn}
Then, $\tilde \Sigma_e$ can be expressed in terms of $\Sigma_e$
as \al{ \tilde \Sigma_e=\left\{
                       \begin{array}{ll}
\Gamma_e\left(I_n-\frac{1-\tau}{\lambda} \Gamma_e^*\Gamma_e\right)^{\frac{1}{\tau-1}}\Gamma_e^*, & \hbox{$0\leq \tau<1$} \\
\Gamma_e\exp\left(\frac{1}{\lambda} \Gamma_e^*\Gamma_e\right)\Gamma_e^*, & \hbox{$\tau=1$,}
                       \end{array}
                     \right.
\nn} where $\Sigma_e=\Gamma_e\Gamma_e^*$. Here, $\lambda$, with $\lambda> (1-\tau)\|\Sigma_e\|_\infty$,
is the unique {\em Lagrange} multiplier
such that $\Sc_\tau(\tilde f\| f)=c$.  \eteo

 Therefore, the noncausal Wiener filter is robust with respect to the $\tau$ divergence constraint. Similarly to the
static case, the worst situation occurs when all the mismodeling budget is allocated in a perturbation of the covariance matrix function $\Sigma_x$.
The connection with the results in \cite{LEVY_NIKOUKHAH_2004} is analogous to the static case.
Moreover, the additional MSE occasioned by the least favorable model perturbation is
$\Delta\mathrm{MSE}=\tr\int (\tilde \Sigma_e-\Sigma_e)$.

Problem (\ref{minimax_wiener}) can be relaxed in the following way 
\al{\label{minimax_dyn_relax}\underset{g\in\Gc}{\min}\;\underset{\tilde f\in\Bc_\tau^\infty}{\max} \;J(\tilde f,g)+\lambda (c-\Sc_\tau(\tilde f\| f)),}
where $\mathcal{B}_\tau^\infty=\{ \tilde f \hbox{ s.t. } \Sc_\tau(\tilde f\| f)<\infty\}.$ Here, $\lambda$ is fixed {\em a priori}
and such that $\lambda>(1-\tau)\|\Sigma_e\|_\infty$. 

\begin{remarkk} Also in this case, it is possible to show that Problem (\ref{minimax_dyn_relax}) is a minimum 
entropy problem in the sense of \cite{MUSTAFA_GLOVER_1990}.
\end{remarkk}
\section{Simulation study}\label{sec_sim}

We analyze the impact of parameters $c$ and $\tau$ on the least favorable statistics of the estimation error
corresponding to the least favorable model in $\Bc_\tau$.
In the static case, given $\tau$ and $c$, $\lambda$ is given by solving equation (\ref{D_tau_lambda})-(\ref{D_1_lambda}) in Appendix. The latter computation can be efficiently performed using the bisection method.    Then, the least favorable statistics of the estimation error is given by Theorem \ref{theorem_static}. The same strategy is applied in the dynamic case. 

\subsection{Static Estimaton}
\begin{figure}[htbp]
\begin{center}
\includegraphics[width=\columnwidth]{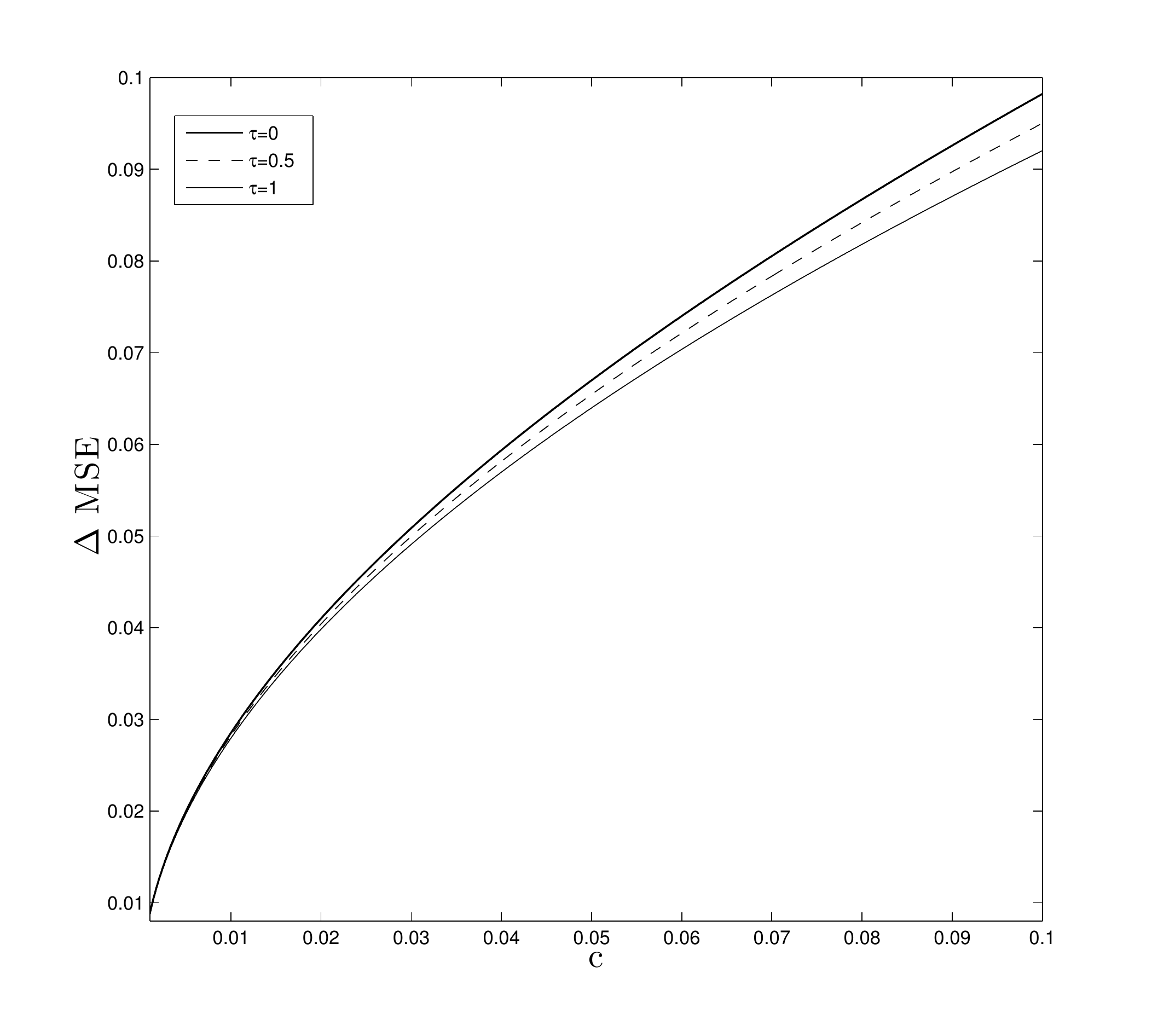}
\end{center}
 \caption{Additional $\mathrm{MSE}$ occasioned by the least favorable perturbation in $\Bc_\tau$ defined in (\ref{def_Bc_static})
 with $\tau=0$, $\tau=0.5$ and $\tau=1$.}\label{dMSE_static}
\end{figure}

We consider a bidimensional Gaussian random vector $x$, that is $n=2$.
We assume that the nominal covariance matrix of the estimation error is 
\al{ P =\left[
     \begin{array}{cc}
      0.15  & 0.05 \\
      0 .05 & 0.1 \\
     \end{array}
   \right].\nn}
 We consider the least favorable statistics
$\tilde f^\circ \in\Bc_\tau$ where $\Bc_\tau$ has been defined in (\ref{def_Bc_static}) with $\tau=0$, $\tau=0.5$ and $\tau=1$. In Figure \ref{dMSE_static} we show the additional MSE occurred
when $c\in[0.001,0.1]$.
We see that, for each value of $\tau$, the larger $c$ is, the larger $\Delta \mathrm{MSE}$ is, as expected.
Moreover, for $c$ fixed, the larger $\tau$ is, the smaller $\Delta \mathrm{MSE}$ is.
Therefore, smaller values of $\tau$ corresponds to least favorable statistics with a larger $\Delta \mathrm{MSE}$.
In order to compare the features of the three different balls $\Bc_\tau$, with $\tau=0$, $\tau=0.5$ and $\tau=1$,
we fix the tolerance $c$ for each ball in such a way that $\Delta \mathrm{MSE}=0.08$ for the least favorable statistics.
The tolerances, respectively, are
$c_{\tau=0}=0.0692$, $c_{\tau=0.5}=0.0728$ and $c_{\tau=1}=0.0767$.
The corresponding covariance matrices of the least favorable estimation error are
\al{ \tilde P_{\tau=0} &=\left[
     \begin{array}{cc}
      0.2041  & 0.0783 \\
      0.0783  &  0.1259 \\
     \end{array}
   \right],\tilde P_{\tau=0.5} =\left[
     \begin{array}{cc}
     0.2039     & 0.0779 \\
      0.0779  & 0.1261 \\
     \end{array}
   \right],\nn\\ \tilde P_{\tau=1} &=\left[
     \begin{array}{cc}
     0.2037   &  0.0775\\
    0 .0775   &  0.1263\\
     \end{array}
   \right].\nn } One can see that the least favorable statistics in $\Bc_{\tau=0}$
   tends to concentrate the perturbation on the component with larger nominal variance. On the contrary,
   the least favorable statistics in $\Bc_{\tau=1}$ tends to spread such perturbation among the two components.
   Finally, the least favorable statistics in $\Bc_{\tau=0.5}$ mitigates those two features.

\subsection{Noncausal Filtering}

\begin{figure}
\begin{center}
\includegraphics[width=\columnwidth]{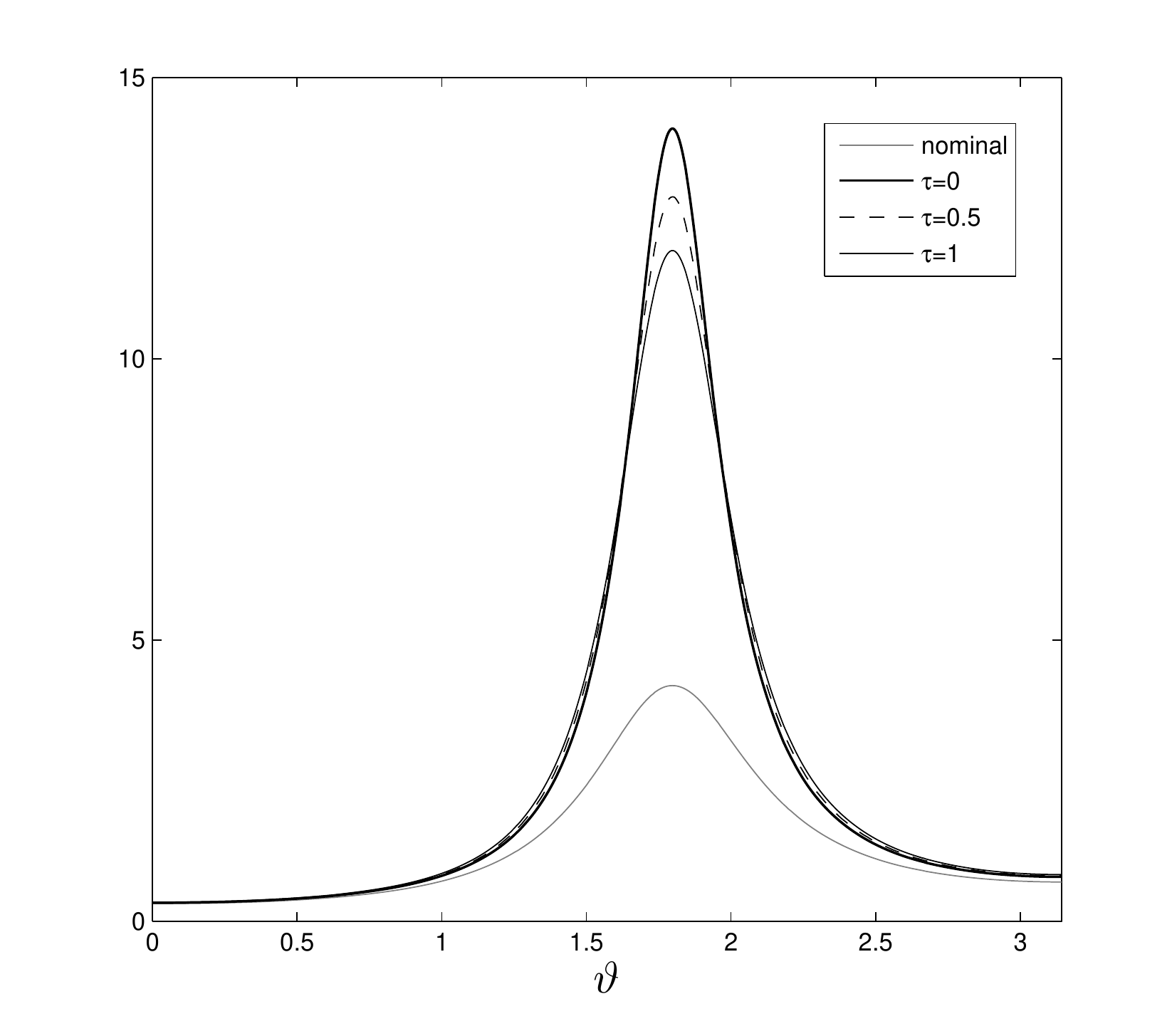}
\end{center}
\caption{Spectral density of $e(t)$ with respect to the nominal statistics and with
respect to the least favorable statistics in $\Bc_\tau$ defined in (\ref{def_Bc_dynamic})
 with $\tau=0$, $\tau=0.5$ and $\tau=1$.}\label{cf_dyn_scalar}
\end{figure}

\begin{figure}
\begin{center}
\includegraphics[width=\columnwidth]{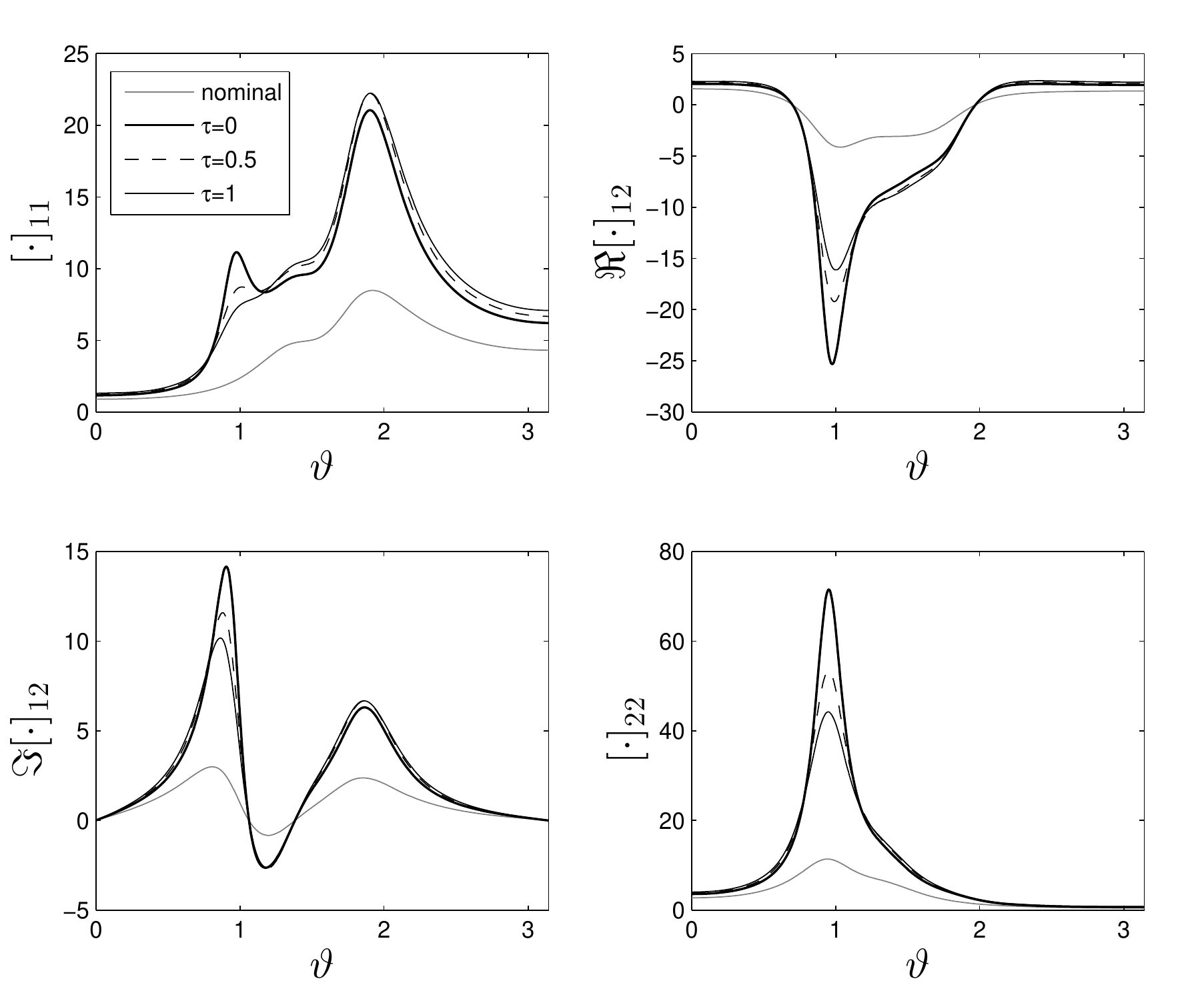}
\end{center}
\caption{Spectral density of $e(t)$ with respect to the nominal statistics and with
respect to the least favorable statistics in $\Bc_\tau$ defined in (\ref{def_Bc_dynamic})
 with $\tau=0$, $\tau=0.5$ and $\tau=1$.}\label{cf_dyn_biv}
\end{figure}

We consider a Gaussian scalar process $x(t)$, that is $n=1$.
The spectral density of its nominal estimation error $e(t)$ is depicted in Figure \ref{cf_dyn_scalar}.
 We consider three different balls $\Bc_\tau$, see (\ref{def_Bc_dynamic}), with $\tau=0$, $\tau=0.5$
and $\tau=1$. Also in this case we noticed that smaller values of $\tau$ corresponds to least favorable probability measures with a larger $\Delta \mathrm{MSE}$.
Similarly to the static case, we fix the tolerance $c$ for each  ball in such a way that $\Delta \mathrm{MSE=0.2}$ for the least favorable statistics.
The tolerances, respectively, are
$c_{\tau=0}=0.022$, $c_{\tau=0.5}=0.025$ and $c_{\tau=1}=0.028$.
The corresponding least favorable spectral densities of $e(t)$ are depicted in Figure
\ref{cf_dyn_scalar}. One can see that the least favorable spectral density in $\Bc_{\tau=0}$
   tends to concentrate the perturbation on the frequency band where the spectral density takes larger values. On the contrary,
   the least favorable spectral density in $\Bc_{\tau=1}$ tends to spread such perturbation over the entire frequency band.
   Also in this case, the least favorable spectral density in $\Bc_{\tau=0.5}$ mitigates those two features. Next, we consider a bidimensional Gaussian process $x(t)$, i.e. $n=2$. The nominal spectral density of $e(t)$
is depicted in Figure \ref{cf_dyn_biv}.
 As before, we fix $c$ for each ball in such a way that the least favorable statistics is such that $\Delta \mathrm{MSE}=1.7$.
We found $c_{\tau=0}=0.076$, $c_{\tau=0.5}=0.088$ and $c_{\tau=1}=0.1$.
The corresponding least favorable
spectral density in $\Bc_\tau$, with $\tau=0$, $\tau=0.5$ and $\tau=1$, is depicted in Figure \ref{cf_dyn_biv}. One can see
that the least favorable spectral density in $\Bc_{\tau=0}$ allocates most of the perturbation in the second component, more precisely
in the frequency band where the nominal spectral density takes the largest values (also compared with respect to the first component). On the contrary, the least favorable spectral density 
in $\Bc_{\tau=1}$ allocates more perturbation in the first component than the former; again, the latter tends to spread the perturbation among the two components and on the entire frequency
band.

\section{Conclusions}\label{sec_concl}
In this paper, we showed that the {\em Bayes} estimator and the noncausal {\em Wiener} filter based on the nominal statistics are robust according to the minimax approach where the uncertainty is specified by a ball formed by placing a bound on the $\tau$ divergence family between the actual and the nominal statistics. 
Interestingly, the relaxation of this minimax problem can be understood as a family of risk-sensitive estimation problems which is also linked 
to a family of minimum entropy problems. 
 Finally, through a simulation study, we have analyzed the features of this family of uncertainty classes. These results represent the starting point to derive a new family of robust {\em Kalman} filters characterized by the $\tau$ divergence family \cite{STATETAU}. 

\section*{Appendix}
 \appendix
 \section{Proof of Proposition \ref{prop_tau_divergence} and Proposition \ref{prop_d_tau_dyn}}
 To prove  Proposition \ref{prop_tau_divergence} it is sufficient to note that \al{\label{tau_div_decomposition}\Dc_\tau(\tilde f\|f)=\frac{1}{1-\tau} \|\Delta m_z\|^2_{K_Z^{-1}}+\overline{\Dc}_\tau(\tilde  K_z\| K_z) ,}
where $\overline{\Dc}_\tau$ is the $\tau$ divergence between $\tilde K_z$ and $K_z$ \cite{OPTIMAL_PREDICTION_ZORZI_2014}.
The latter is nonnegative and is equal to zero if and only if $\tilde K_z=K_z$. Since $K_z\in\Qc_+^{q}$, the first term in (\ref{tau_div_decomposition}) is nonnegative and  is equal to zero if and only if $\Delta m_z=0$, that is $\tilde m_z=m_z$. Therefore, $\Dc_\tau$ is nonnegative and equality holds if and only if $\tilde m_z=m_z$ and $\tilde K_z=K_z$, that is $\tilde f=f$.
The proof for cases $\tau=0$ and $\tau=1$ is similar. Finally, Proposition  \ref{prop_d_tau_dyn} can be proved along the same lines.  \qed

\section{Proof of Theorem \ref{theorem_static}}
The proof is divided in three cases.

\textbf{Case $\mathbf{\tau=0}$.} $\Dc_\tau$ is equivalent to the {\em Kullback-Leibler} divergence and the statement has been proved in \cite[Theorem 1]{LEVY_NIKOUKHAH_2004}.\\

 \textbf{Case $\mathbf{0<\tau <1}$.} We have to show that  $ J(\tilde f,g^\circ )\leq J(\tilde f^\circ,g^\circ)\leq J(\tilde f^\circ,g)$
for any $(\tilde f,g)\in \Bc_\tau \times \Gc$. Since $\tilde f^\circ\sim \Nc(\tilde m_z^\circ,\tilde K_z^\circ)$,
the inequality $J(\tilde f^\circ,g^\circ)\leq J(\tilde f^\circ,g)$ implies that $g^\circ$
is the Bayesian estimator (\ref{Opt_estimator}). Next, we show the inequality $J(\tilde f,g^\circ )\leq J(\tilde f^\circ,g^\circ)$ holds. Therefore, it is sufficient to show
 \al{  \tilde f^\circ=\underset{\tilde f\in\Bc_{\tau}}{\arg \max}\, J(\tilde f,g^\circ).}
 Let $ e=x-g^\circ(y)=\left[\begin{array}{cc} I_n & -G^\circ \\ \end{array}\right](z-m_z)$.
 Therefore, $ \tilde m_e:=\Es_{\tilde f}[e]=\left[\begin{array}{cc} I_n & -G^\circ \\ \end{array}\right]\Delta m_z$
and \al{ \tilde K_e&:= \Es_{\tilde f}[(e-\tilde m_e)(e-\tilde m_e)^T]\nn\\ &=\left(\left[\begin{array}{cc} I_n & -G^\circ \\ \end{array}\right] \tilde K_z \left[\begin{array}{c}I_n \\ -(G^\circ)^T\end{array}\right]\right) .\nn }
 Moreover,
\al{ J( &\tilde f,g^\circ)=\tr (\tilde  m_e\tilde m_e^T+ \tilde K_e)\nn\\
&= \tr\left( \left[
                             \begin{array}{cc}
                               I_n & -G^\circ \\
                             \end{array}
                           \right] (\tilde K_z+\Delta m_z \Delta m_z^T)\left[
                                                    \begin{array}{c}
                                                      I_n\\
                                                      -(G^\circ)^T \\
                                                    \end{array}
                                                  \right]
\right).\nn}
In order to characterize $\tilde f$, we exploit the duality theory. The Lagrangian is
\al{ \label{lagrangian_static} \Lc &(\tilde m_z,\tilde K_z,\lambda)=J(\tilde f,g^\circ )+\lambda(c-\Dc_\tau(\tilde f\| f))\nn\\
& =\Delta m_z^T W_\lambda \Delta m_z +\tr\left(\left[
                                                    \begin{array}{cc}
                                                      I_n & -G^\circ \\
                                                    \end{array}
                                                  \right] \tilde K_z
\left[
                                                    \begin{array}{c}
                                                       I_n \\
                                                       -(G^\circ)^T \\
                                                    \end{array}
                                                  \right]
\right)\nn\\ &\hspace{0.3cm}+\lambda\left(c+\tr\left(\frac{1}{\tau(1-\tau)}(L_z^{-1} \tilde{K}_z L_z^{-T})^{\tau}\right.\right.\nn\\
& \hspace{0.3cm} \left. \left.-\frac{1}{1-\tau} \tilde{K}_z K_z^{-1} -\frac{1}{\tau} I_{n+p} \right)\right),}
where
\al{ W_\lambda=
\left[
                                                    \begin{array}{c}
                                                       I_n \\
                                                       -(G^\circ)^T \\
                                                    \end{array}
                                                  \right] \left[
                                                    \begin{array}{cc}
                                                      I_n & -G^\circ \\
                                                    \end{array}
                                                  \right]  -\frac{\lambda}{1-\tau} K_z^{-1}\nn }
                                                  and $\lambda\geq 0$ is the {\em Lagrange} multiplier.
Note that, $\Lc$ is bounded above and strictly concave in $\tilde K_z$ when $\lambda>0$. Moreover, $\Lc$ is bounded above and strictly concave in $\tilde m_z$ if $W_\lambda$ is negative definite.
Define $ M=\left[
            \begin{array}{cc}
             I_n & -G^\circ \\
            \end{array}
          \right]L_z$. Since $G^\circ=K_{xy}K_y^{-1}$ and  $K_z=L_zL_z^T$, it is easy to check that
\al{ MM^T&=\left[
                                                    \begin{array}{cc}
                                                     I_n & -G^\circ\\
                                                    \end{array}
                                                  \right]K_z\left[
                                                    \begin{array}{c}
                                                       I_n \\
                                                       -(G^\circ)^T \\
                                                    \end{array}
                                                  \right]\nn\\
                                                  &= K_{x}-K_{xy} K_y^{-1}K_{yx}=P. \nn}
Moreover,  $W_\lambda$ is congruent to $ M^TM-\frac{\lambda}{1-\tau} I_{n+p}$
which is negative definite when \al{\label{cond_lambda} \lambda> (1-\tau)\|M^T M\|=(1-\tau) \|MM^T\|=(1-\tau)\|P\|.}
Therefore, under assumption (\ref{cond_lambda}),  $\Lc$ is bounded above  and strictly concave in  $(\tilde  m_z,\tilde K_z)$, and it is maximized by the points $(\tilde m_z^\circ,\tilde K_z^\circ)$ annihilating its first variations
\al{ \nabla _{\tilde m_z,u}\Lc &=2\Delta m_z^T W_\lambda u\nn\\  \nabla_{\tilde K_z,V} \Lc &=\tr\left(\left( \left[
                                                    \begin{array}{c}
                                                        I_n \\
                                                        -(G^\circ)^T\\
                                                    \end{array}
                                                  \right]\left[
                                                    \begin{array}{cc}
                                                      I_n  & -G^\circ\\
                                                    \end{array}
                                                  \right]\right. \right. \nn\\   & \hspace{-0.3cm}\left.\left.                                                 +\frac{\lambda}{1-\tau}(L_z^{-T}(L_z^{-1}\tilde K_z L_z^{-T})^{\tau-1}L_z^{-1}-K_z^{-1})\right) V\right)\nn}
for any direction $u\in\Rs^{n+p}$ and $V\in\Qc^{n+p}$, respectively. To compute $\nabla_{\tilde K_z,V}\Lc$ we exploited the formula for the first variation 
of the exponentiation  of a matrix given in \cite{BETA}.
Since $W_\lambda$ is negative definite, it follows that $\Delta m_z=0$, and thus $\tilde m_z^\circ=m_z$. Regarding $\tilde K_z^\circ$, we have
\al{ \frac{\lambda}{1-\tau}  L_z^{-T} & (L_z^{-1}\tilde K_z^\circ L_z^{-T})^{\tau-1}L_z^{-1}\nn\\ &=\frac{\lambda}{1-\tau} K_z^{-1}-\left[
                                                    \begin{array}{c}
                                                      I_n\\
                                                      -(G^\circ)^T  \\
                                                    \end{array}
                                                  \right]\left[
                                                    \begin{array}{cc}
                                                       I_n & -G^\circ \\
                                                    \end{array}
                                                  \right] \nn} moreover,
\al{   \label{forma_parziale_ott}(L_z^{-1} & \tilde K_z^\circ L_z^{-T})^{\tau-1}\nn \\ & =I_{n+p}-\frac{1-\tau}{\lambda }L_z^T\left[
                                                    \begin{array}{c}
                                                      I_n \\
                                                     -(G^\circ)^T \\
                                                    \end{array}
                                                  \right]\left[
                                                    \begin{array}{cc}
                                                       I_n & -G^\circ  \\
                                                    \end{array}
                                                  \right]L_z,}  where the right hand side of (\ref{forma_parziale_ott}) is positive definite by (\ref{cond_lambda}).
Accordingly, under assumption (\ref{cond_lambda}), we have
\al{\label{Kz_tilde_optimal} & \tilde K_z^\circ  =\nn\\ & \hspace{0.0cm}L_z\left(I_{n+p}-\frac{1-\tau}{\lambda} L_z^T\left[
                                                      \begin{array}{c}
                                                        I_n\\
                                                          -(G^\circ)^T  \\
                                                      \end{array}
                                                    \right]\left[
                                                             \begin{array}{cc}
                                                               I_n &  -G^\circ \\
                                                             \end{array}
                                                           \right]L_z
\right)^{\frac{1}{\tau-1}} L_z^T.} Note that, $K_z$ admits the following block upper diagonal lower (UDL) factorization
\al{ \label{fact_Kz} K_z=\left[
       \begin{array}{cc}
         I_{n} & G^\circ \\
         0 & I_{p} \\
       \end{array}
     \right]\left[
       \begin{array}{cc}
         P & 0 \\
         0 & K_y \\
       \end{array}
     \right]\left[
       \begin{array}{cc}
         I_{n} & 0 \\
         (G^\circ)^T & I_{p} \\
       \end{array}
     \right]}
and we choose $L_z$ as \al{ \label{def_Lz}L_z=\left[
              \begin{array}{cc}
                I_{n} & G^\circ \\
               0  & I_{p} \\
              \end{array}
            \right]\left[
                     \begin{array}{cc}
                       L_P & 0 \\
                       0 & L_y \\
                     \end{array}
                   \right],
} where $K_y=L_yL_y^T$ and $P=L_P L_P^T$.
Accordingly,
\al{  & \left(I_{n+p}-\frac{1-\tau}{\lambda} L_z^T \left[
                                   \begin{array}{c}
                                     I_n \\
                                     -(G^\circ)^T \\
                                   \end{array}
                                 \right]\left[
                                          \begin{array}{cc}
                                            I_n & -G^\circ \\
                                          \end{array}
                                        \right]
L_z\right)^{\frac{1}{\tau-1}}\nn\\ &\hspace{0.2cm}=\left(I_{n+p}-\frac{1-\tau}{\lambda}\left[
      \begin{array}{c}
        L_P^T \\
       0  \\
      \end{array}
    \right]\left[
             \begin{array}{cc}
               L_P & 0  \\
             \end{array}
           \right]\right)^{\frac{1}{\tau-1}}\nn\\ & \hspace{0.2cm}=\left(\left[
                                         \begin{array}{cc}
                                           I_{n}-\frac{1-\tau}{\lambda} L_P^T L_P  & 0 \\
                                           0 & I_p \\
                                         \end{array}
                                       \right]\right)^{\frac{1}{\tau-1}}\nn\\
& \label{manipulations1}\hspace{0.2cm}=\left[
                                         \begin{array}{cc}
                                           (I_{n}-\frac{1-\tau}{\lambda} L_P^T L_P)^{\frac{1}{\tau-1}} & 0 \\
                                           0 & I_p \\
                                         \end{array}
                                       \right]}
and by (\ref{Kz_tilde_optimal}), (\ref{def_Lz}) and (\ref{manipulations1}) we have
\al{ \label{fact_tilde_Kz} \tilde K_z^\circ=\left[
       \begin{array}{cc}
         I_{n} & G^\circ \\
         0 & I_{p} \\
       \end{array}
     \right]\left[
       \begin{array}{cc}
         \tilde P & 0 \\
         0 & K_y \\
       \end{array}
     \right]\left[
       \begin{array}{cc}
         I_{n} & 0 \\
         (G^\circ)^T & I_{p} \\
       \end{array}
     \right],}
     where \al{ \label{Ptilda} \tilde P=L_P\left(I_{n}-\frac{1-\tau}{\lambda} L_P^T L_P\right)^{\frac{1}{\tau-1}} L_P^T.}
     Moreover,
\al{ \tilde K_z^\circ=\left[
                 \begin{array}{cc}
                   \tilde P+K_{xy}K_y^{-1}K_{yx} & K_{xy} \\
                   K_{yx} &  K_y \\
                 \end{array}
               \right]
\nn}
accordingly, $  \tilde K_x=\tilde P+K_{xy}K_y^{-1}K_{yx} $.
Since $\lambda>0$, the duality gap between the primal and the dual is zero if and only if $\Dc_\tau(\tilde f^\circ\| f)=c$. It remains to be shown that there exists $\lambda>(1-\tau)\|P\|$ such that $\Dc_\tau(\tilde f^\circ\| f)=c$.
By considering the factorizations (\ref{fact_Kz}), (\ref{fact_tilde_Kz}) and (\ref{Ptilda}), we obtain
\al{ \label{D_tau_lambda}& \Dc_\tau(\tilde f^\circ\| f)\nn\\ &=\tr\left(\frac{1}{\tau(\tau-1)}(L_z^{-1}\tilde K_z^\circ L_z^{-T})^{\tau}+\frac{1}{1-\tau}\tilde K_z^\circ K_z^{-1}+\frac{1}{\tau}I_{n+p}\right)\nn\\
& = \tr\left(\frac{1}{\tau(\tau-1)}(L_P^{-1}\tilde P L_P^{-T})^{\tau}+\frac{1}{1-\tau} L_P^{-1}\tilde P L_P^{-T}+\frac{1}{\tau}I_{n}\right)\nn\\
 & = \tr\left(\frac{1}{\tau(\tau-1)}\left(I_{n}-\frac{1-\tau}{\lambda}L_P^T L_P\right)^{\frac{\tau}{\tau-1}}\right.\nn\\ & \hspace{0.3cm}\left. +\frac{1}{1-\tau}\left(I_{n}-\frac{1-\tau}{\lambda}L_P^T L_P\right)^{\frac{1}{\tau-1}}+\frac{1}{\tau}I_{n}\right).}
Let $L_P^TL_P=UDU^T$ be the eigenvalue decomposition of $L_P^TL_P$, where $UU^T=I_n$ and
$D=\mathrm{diag}(d_1,d_2,\ldots,d_n)$ with $d_i=\sigma_i(L_P^T L_P)$.
Therefore, we get
\al{  \Dc_\tau(\tilde f^\circ & \| f)\nn\\ & = \tr\left(\frac{1}{\tau(\tau-1)}\left(I_{n}-\frac{1-\tau}{\lambda}D\right)^{\frac{\tau}{\tau-1}}\right. \nn\\& \hspace{0.3cm} \left. +\frac{1}{1-\tau}\left(I_{n}-\frac{1-\tau}{\lambda}D\right)^{\frac{1}{\tau-1}}+\frac{1}{\tau}I_{n}\right)\nn\\
&= \sum_{i=1}^n \gamma(\lambda,d_i)
\nn }
with \al{ \gamma(\lambda,d_i)=&\frac{1}{\tau(\tau-1)} \left(1-\frac{1-\tau}{\lambda}d_i\right)^{\frac{\tau}{\tau-1} \nn}
\nn\\ &  +\frac{1}{1-\tau}\left(1-\frac{1-\tau}{\lambda}d_i\right)^{\frac{1}{\tau-1}}+\frac{1}{\tau}. \nn }
For $\lambda>(1-\tau)d_i>0$, we have
\al{\frac{d}{d\lambda}\gamma(\lambda,d_i)=-\frac{d_i^2}{\lambda^3}\left(1-\frac{1-\tau}{\lambda} d_i\right)^{\frac{2-\tau}{\tau-1}}<0 \nn }
and \al{ \underset{\lambda\rightarrow \infty}{\lim}\gamma(\lambda,d_i)=0,\;\; \underset{\lambda\rightarrow ((1-\tau)d_i)^+}{\lim}\gamma(\lambda,d_i)=\infty \nn }
so that $\gamma(\lambda,d_i)$ is a monotone decreasing function of $\lambda$. Since $\Dc_{\tau}(\tilde f^\circ\| f)$ is the the sum of $\gamma(\lambda,d_i)$ with $i=1\ldots n$, it is also a monotone decreasing
function of $\lambda$ over $((1-\tau)\|P\|,\infty)$ and
\al{ \underset{\lambda\rightarrow \infty}{\lim}\Dc_{\tau}(\tilde f^\circ\| f)=0,\;\; \underset{\lambda\rightarrow ((1-\tau)\|P\|)^+}{\lim}\Dc_{\tau}(\tilde f^\circ\| f)=\infty.} We conclude, for any $c>0$ there exists a unique $\lambda>(1-\tau)\|P\|$ such that $\Dc_\tau(\tilde f^\circ\| f)=c$.\\

\textbf{Case $\mathbf{\tau=1}$.} The Lagrangian $\Lc$ can be formed as in (\ref{lagrangian_static}). Then, one can see that $\Lc$ is bounded below if and only if $\tilde m_z= m_z$. Accordingly, the least favorable density $\tilde f^\circ$ has mean $m_z$. Then, the least favorable $\tilde K_z^\circ$ can be characterized using the duality theory similarly  to the case $0<\tau <1$. In particular, it is not difficult to see that
\al{ \label{D_1_lambda}\Dc_1(& \tilde f^\circ\| f)\nn\\ &=\tr\left(\exp\left(\frac{1}{\lambda} L_P^T L_P\right)\left(\frac{1}{\lambda}L_P^T L_P-I_n \right) +I_n\right). }
\qed

\section{Proof of Proposition \ref{prop_tau_entropy}}
The statement can be proved along the same lines of Proposition 2.3.1 and Proposition 2.3.2 in \cite{MUSTAFA_GLOVER_1990}.\qed

\section{Proof of Theorem \ref{TH_min_entropy}}
Exploiting similar argumentations used in the proof of Proposition 1 in \cite{LEVY_NIKOUKHAH_2004}, we have that
 \al{ \max_{\tilde f\in\Bc_\tau^\infty} \Lc(\tilde f,\lambda,g)=\Hc_\tau(e,\lambda)+\lambda c,  \nn} 
where  $\Lc(\tilde f,\lambda,g)=J(\tilde f,g)+\lambda(c-\Dc_\tau(\tilde f\| f ))$. Since $g$ does not depend on $\lambda$, we conclude that 
\al{\min_{g\in\Gc} \max_{\tilde f\in\Bc_\tau^\infty} \Lc(\tilde f,\lambda,g)=\min_{g\in\Gc}\Hc_\tau(e,\lambda). \hspace{0.6cm}\qed \nn }

\section{Proof of Theorem \ref{thm_wiener_saddle_point}}
The proof is similar to the one of Theorem \ref{theorem_static}. The unique peculiarity follows. 
Since  $\Sigma_z$ has bounded {\em McMillan} degree and $\frac{1}{1-\tau}\in\Ns$, then the integrand function in $\Sc_{\tau}(\tilde f^\circ\| f)$ is rational. Moreover,
 as $\lambda\rightarrow ((1-\tau)\|\Sigma_e\|_\infty)^+$
 this integrand function has at least one pole tending to $\Ts$. Accordingly, these assumptions allow to conclude that 
\al{ \underset{\lambda\rightarrow ((1-\tau)\|\Sigma_e\|_\infty)^+}{\lim}\Sc_{\tau}(\tilde f^\circ\| f)=\infty.\hspace{0.6cm}\qed	\nn }


\end{document}